\documentclass[12pt,a4paper,twoside,reqno]{amsart}

\usepackage{amsmath, amssymb, amsthm, enumerate, enumitem, graphicx}
\usepackage[normalem]{ulem}
\usepackage{xr}
\makeatletter
\def\l@section{\@tocline{1}{10pt}{1pc}{}{}}
\def\l@subsection{\@tocline{2}{0pt}{1pc}{4.6em}{}}
\def\l@subsubsection{\@tocline{3}{0pt}{1pc}{7.6em}{}}
\renewcommand{\tocsection}[3]{%
  \indentlabel{\@ifnotempty{#2}{\makebox[2.3em][l]{%
    \ignorespaces#1 #2.\hfill}}}\textbf{#3}}
\renewcommand{\tocsubsection}[3]{%
  \indentlabel{\@ifnotempty{#2}{\hspace*{2.3em}\makebox[2.3em][l]{%
    \ignorespaces#1 #2.\hfill}}}#3}
\renewcommand{\tocsubsubsection}[3]{%
  \indentlabel{\@ifnotempty{#2}{\hspace*{4.6em}\makebox[3em][l]{%
    \ignorespaces#1 #2.\hfill}}}#3}
\makeatother

\newcommand{\MM}{\mathcal{M}}

\newcommand{\IR}{\mathbb{R}}

\newcommand{\IH}{\mathbb{H}}

\newcommand{\IZ}{\mathbb{Z}}

\newcommand{\TT}{\mathcal{T}}

\newcommand{\ov}[1]{\overline{#1}}

\newcommand{\td}[1]{\widetilde{#1}}
\DeclareMathOperator{\Nil}{Nil}
\DeclareMathOperator{\Klein}{Klein}

\DeclareMathOperator{\stan}{stan}

\DeclareMathOperator{\Int}{Int}

\DeclareMathOperator{\Ric}{Ric}

\DeclareMathOperator{\area}{area}

\DeclareMathOperator{\diam}{diam}

\DeclareMathOperator{\vol}{vol}

\DeclareMathOperator{\Rm}{Rm}

\newcommand{\dotcup}{\ensuremath{\mathaccent\cdot\cup}}
\newcommand{\EMPTY}[1]{}

\newtheorem{Theorem}{Theorem}[section]

\newtheorem{Corollary}[Theorem]{Corollary}

\newtheorem{Question}[Theorem]{Question}
\numberwithin{equation}{section}

\newtheorem*{KeyLemma}{Key Lemma}
\newtheorem*{GeneralizedKeyLemma}{Generalized Key Lemma}
\theoremstyle{definition} 
\newtheorem*{paperA}{A: Generalizations of Perelman's long-time estimates ([BamA])}
\newtheorem*{paperB}{B: Evolution of the minimal area of simplicial complexes under Ricci flow ([BamB])}
\newtheorem*{paperC}{C: 3-manifold topology and combinatorics of simplicial complexes in $3$-manifolds ([BamC])}
\newtheorem*{paperD}{D: Proof of the main results ([BamD])}

\title[Long-time behavior of 3d Ricci flow --- Introduction]{Long-time behavior of 3 dimensional Ricci flow\\Introduction}
\author{Richard H Bamler}
\address{UC Berkeley, Department of Mathematics, 970 Evans Hall, Berkeley, CA 94720, USA}
\email{rbamler@math.berkeley.edu}
\date{\today}

\begin{document}
\begin{abstract}
In the following series of papers we analyze the long-time behavior of 3 dimensional Ricci flows with surgery.
Our main result will be that if the surgeries are performed correctly, then only finitely many surgeries occur and after some time the curvature is bounded by $C t^{-1}$.
This result confirms a conjecture of Perelman.
In the course of the proof, we also obtain a qualitative description of the geometry as $t \to \infty$.
\end{abstract}

\maketitle
\tableofcontents

\section{Introduction} 
\subsection{Statement of the main result}
In the following series of papers, we analyze the long-time behavior of $3$ dimensional Ricci flows with surgery and we prove a conjecture of Perelman.
In a few words, our first main result can be summarized as follows.
We refer to Theorem \ref{Thm:MainTheorem-III} on page \pageref{Thm:MainTheorem-III} for a precise statement.

\begin{quote}
\textit{Let $(M,g)$ be a closed and orientable $3$ dimensional Riemannian manifold. \\
Then there is a Ricci flow with only  \emph{finitely} many surgeries whose initial time-slice is $(M, g)$ and that either goes extinct in finite time or exists for all times.
Moreover, there is a constant $C$ such that the norm of the Riemannian curvature tensor in this flow is bounded everywhere by $C t^{-1}$ for large times $t$.}
\end{quote}
We moreover obtain a characterization of the geometry of this Ricci flow at large times, which will be summarized in Theorem \ref{Thm:geombehavior} on page \pageref{Thm:geombehavior}.

The Ricci flow with surgery has been used by Perelman to solve the Poincar\'e and Geometrization Conjectures (\cite{PerelmanI}, \cite{PerelmanII}, \cite{PerelmanIII}).
Given any initial metric on a closed $3$-manifold, Perelman managed to construct a solution to the Ricci flow with surgery on a maximal time-interval and showed that its surgery times do not accumulate.
This means that every finite time-interval contains only a finite number of surgery times.
Furthermore, he proved that if the given manifold is a homotopy sphere (or, more generally, a connected sum of prime, non-aspherical manifolds), then this flow goes extinct in finite time and the total number of surgeries is finite.
This fact implies that the initial manifold is a sphere if it is simply connected and hence proves the Poincar\'e Conjecture.
On the other hand, if the Ricci flow with surgery continues to exist for infinite time, Perelman showed that the manifold decomposes into a thick part, which approaches a hyperbolic metric, and a thin part, which becomes arbitrarily collapsed on local scales.
Based on this collapse, it is then possible to show that the thin part can be decomposed into pieces whose topology is well understood (\cite{ShioyaYamaguchi}, \cite{MorganTian}, \cite{KLcollapse}).
Eventually, this decomposition can be reorganized to a geometric decomposition, thus establishing the Geometrization Conjecture.

Observe that although the Ricci flow with surgery was used to solve such hard problems, some of its basic properties are still unknown, because they surprisingly turned out to be irrelevant in the end.
For example, it was only conjectured by Perelman that in the long-time existent case there are finitely many surgeries, i.e. that after some time the flow can be continued by a conventional smooth, non-singular Ricci flow defined up to time infinity.
Furthermore, it is still unknown whether and in what way precisely the Ricci flow exhibits the full geometric decomposition of the manifold.

In \cite{LottTypeIII}, \cite{LottDimRed} and \cite{LottSesum}, Lott and Lott-Sesum gave a description of the long-time behavior of certain non-singular Ricci flows on manifolds whose geometric decomposition consists of a single component.
However, they needed to make additional curvature and diameter or symmetry assumptions.
In \cite{Bamler-longtime-I}, the author proved that under a purely topological condition (sometimes referred to as $\TT_1$), which roughly states that the manifold only consists of  hyperbolic components, the number of surgeries is indeed finite and the curvature is bounded by $C t^{-1}$ for large $t$.
In this paper we remove this additional condition and only assume that the initial manifold is closed and orientable.

This series of papers is a restructured version of the two preprints \cite{Bamler-longtime-II} and \cite{Bamler-longtime-III}.
In \cite{Bamler-longtime-II}, the condition $\TT_1$ was generalized to a far more general topological condition $\TT_2$, which requires that the non-hyperbolic pieces in the geometric decomposition of the underlying manifold contain sufficiently many incompressible surfaces.
For example, manifolds of the form $\Sigma \times S^1$ for closed, orientable surfaces $\Sigma$, in particular the $3$-torus $T^3$, satisfy property $\TT_2$, but the Heisenberg manifold does not.
We refer to \cite[subsection 1.2]{Bamler-longtime-II} for a precise definition and discussion of the conditions $\TT_1$ and $\TT_2$.
Eventually, in \cite{Bamler-longtime-III} the result was further generalized and condition $\TT_2$ was removed.
This generalization was obtained by replacing said incompressible surfaces by simplicial complexes in a careful way.
In the present paper we have merged the proofs of \cite{Bamler-longtime-II} and \cite{Bamler-longtime-III}, so the conditions $\TT_1, \TT_2$ as well as the incompressible surfaces don't play a role anymore.

We now state our main result.
The notions relating to ``Ricci flows with surgery'' that are used in the following are explained in subsection \ref{sec:DefRFsurg} of the following paper \cite{Bamler-LT-Perelman}.
A Ricci flow $(g_t)_{t \in I}$ on a manifold $M$ is a smooth family of Riemannian metrics that satisfy the evolution equation
\[ \partial_t g_t = - 2 \Ric_{g_t}. \]
A ``Ricci flow with surgery $\MM$, that is performed by $\delta(t)$-precise cutoff'' can briefly be described as a sequence of 3 dimensional Ricci flows $(M^1, (g^1_t)_{t \in [0,T^1)}), \linebreak[1] (M^2, (g^2_t)_{t \in [T^1, T^2)}), \ldots$ such that the time-$T^i$ slice $(M^{i+1}, g^{i+1}_{T^i})$ is obtained from the singular metric $g^i_{T^i}$ on $M^i$ by a so called surgery process, which amounts to a geometric version of an inverse connected sum decomposition at a scale less than $\delta(T^i)$ and the removal of spherical or $S^2 \times S^1$ components.
We allow the case in which there are only finitely many surgery times $T^i$ and in which $T^i = \infty$ for the final index $i$.
Observe that we have chosen our notion such that a  $\delta(t)$-precise cutoff is also $\delta'(t)$-precise if $\delta'(t) \geq \delta(t)$.

In \cite{PerelmanII} Perelman showed the existence of a (non-explicit) function $\delta(t)$ such that every normalized Riemannian manifold $(M,g)$ (see Definition \ref{Def:normalized} in the next paper \cite{Bamler-LT-Perelman}) can be evolved into a Ricci flow with surgery $\MM$ that is performed by $\delta(t)$-precise cutoff\footnote{Perelman uses a slightly different notion of Ricci flow with $\delta(t)$-cutoff. For example, he performs surgeries at a scale $= r(t) \delta^2 (t)$ instead of $< \delta(t)$, where $r(t)$ is another function that goes to $0$ as $t \to \infty$. Both notions are however equivalent modulo the choice of a different function $\delta(t)$.} and he proved that for any such Ricci flow with surgery---performed by $\delta(t)$-precise cutoff and with normalized initial conditions---the sur\-gery times $T^i$ do not accumulate.
So if there were infinitely many surgery times (or, equivalently, infinitely many surgeries), then we must have $\lim_{i \to \infty} T^i = \infty$.
Our main result now states that this cannot happen under normalized initial conditions and if $\delta(t)$ has been chosen sufficiently small.
Note that these two conditions are not very restrictive since they already had been imposed in Perelman's work.

\begin{Theorem} \label{Thm:LT0-main-1}
Given a surgery model $(M_{\stan}, g_{\stan}, D_{\stan})$, there is a continuous function $\delta : [0, \infty) \to (0, 1)$ such that the following holds: \nopagebreak  \label{Thm:MainTheorem-III}

Let $\MM$ be a (3 dimensional) Ricci flow with surgery, defined on some time-interval $[0,T_0)$ for $T_0 \leq \infty$, that has normalized initial conditions and is performed by $\delta(t)$-precise cutoff.

Then $\MM$ has only finitely many surgeries and there are constants $T, C < \infty$ such that $|{\Rm_t}| < C t^{-1}$ on $\MM(t)$ for all $t \geq T$.
\end{Theorem}

Note that this curvature bound is optimal apart from the non-explicit constant $C$.
For example, if we consider a Ricci flow on a hyperbolic manifold, then the sectional curvature behaves like $- \frac14 t^{-1}$ as $t \to \infty$.

We mention two interesting direct consequences of Theorem \ref{Thm:MainTheorem-III}, which can be expressed in a completely elementary way and which illustrate the power of this theorem.
None of these results have been proven so far to the author's knowledge.
The first consequence is just a restatement of Theorem \ref{Thm:MainTheorem-III} in the case in which $\MM$ is non-singular.
Note that even in this particular case our proof does not simplify significantly apart from the fact that we don't have to deal with various technicalities.
In fact, the reader is advised to only consider non-singular Ricci flows upon first reading of this series of papers.

\begin{Corollary} \label{Cor:nonsingflow}
Let $(M, (g_t)_{t \in [0, \infty)})$ be a non-singular, long-time existent Ricci flow on a closed $3$-manifold $M$.
Then there is a constant $C < \infty$ such that
\[ |{\Rm_t}| < \frac{C}{t+1} \qquad \text{for all} \qquad t \geq 0. \]
\end{Corollary}

The next result provides a characterization of when the condition of the previous corollary can indeed be satisfied.
\begin{Corollary} \label{Cor:topconditionforlongtime}
Let $M$ be a closed $3$-manifold.
Then there exists a long-time existent Ricci flow $(g_t)_{t \in [0, \infty)}$ on $M$ if and only if $\pi_2(M) = \pi_3(M) = 0$.
\end{Corollary}

Note that this topological condition is equivalent to $M$ being aspherical, which is equivalent to $M$ being irreducible and not diffeomorphic to a spherical space form.

This corollary can be deduced from Theorem \ref{Thm:MainTheorem-III} as follows:
Any normalized Riemannian metric $g$ on an aspherical manifold $M$, can be evolved into a long-time existent Ricci flow with surgery $\MM$ on the time-interval $[0, \infty)$ that is performed by $\delta(t)$-precise cutoff, due to Perelman (\cite{PerelmanII}, see also \cite[Proposition \ref{Prop:RFwsurg-existence}]{Bamler-LT-Perelman}).
The topological condition ensures that all surgeries on $\MM$ are trivial and hence that every time-slice of $\MM$ has a component that is diffeomorphic to $M$.
By Theorem \ref{Thm:MainTheorem-III}, there is a final surgery time $T < \infty$ on $\MM$.
So the flow $\MM$ restricted to the time-interval $[T, \infty)$ is non-singular and the underlying manifold is diffeomorphic to $M$.
Shifting this flow in time by $-T$ yields the desired Ricci flow.
The reverse direction is well known, for example it is a direct consequence of Proposition \ref{Prop:irreducibleafterfinitetime} in the last paper of series, \cite{Bamler-LT-main}, and finite-time extinction (see \cite{PerelmanIII}, \cite{ColdingMinicozziextinction}, \cite{MTRicciflow}).

In the course of the proof of Theorem \ref{Thm:MainTheorem-III} we will obtain a more detailed description of the geometry of the time-slices $\MM(t)$ for large times $t$.
In short, we will find that as $t \to \infty$ the Ricci flow decomposes the manifold into regions that are hyperbolic or exhibit different collapsing behavior at scale $\sqrt{t}$.
The collapse can either be observed on the whole manifold, in which case we speak of a \emph{total collapse} and the underlying manifold is a quotient of a torus or a Heisenberg manifold.
Or it occurs along incompressible (i.e. $\pi_1$-injective) circle ($S^1$) or $2$-torus ($T^2$) fibers.
A collapse along $S^1$-fibers gives rise to Seifert structures (compare with the white regions $M_{1,t}, \ldots, M_{6,t}$ in Figure \ref{fig:longtimepicture}). 
Regions that collapse along $T^2$-fibers look like $T^2 \times I$ or like a twisted interval bundle over the Klein bottle (compare with the gray regions $E_{1,t}, \ldots, E_{7,t}$ in Figure \ref{fig:longtimepicture}).
Those regions either cover the whole manifold, in which case the manifold is a quotient of a $2$-torus bundle over a circle, or they serve as interpolations between different Seifert fibrations.
By this we mean the following: 
A region $E_{i,t} \approx T^2 \times I$ that collapses along the $T^2$ factor is adjacent to two regions $M_{j_1,t}, M_{j_2,t}$ that carry Seifert fibrations.
The Seifert fibers on $M_{j_1,t}, M_{j_2,t}$ represent possibly $S^1$-different directions the boundary tori of $E_{i,t}$.
Towards the ends of $E_{i,t}$ one $S^1$-direction of the $2$-tori along which $E_{i,t}$ collapses becomes so large that on $M_{j_1,t}$ or $M_{j_2,t}$ we only observe a collapse along the other $S^1$-direction.

This decomposition of the underlying manifold into regions that are hyperbolic or Seifert corresponds to a geometric decomposition as defined in Definition \ref{Def:geomdec} of the third paper in this series, \cite{Bamler-LT-topology}.
However, this decomposition is not necessarily minimal, i.e. it may a priori be possible to simplify it by fusing together certain Seifert structures.

We will now summarize our findings more precisely:

\begin{figure}[t] 
\begin{center}
\setlength{\unitlength}{2863sp}%
\begingroup\makeatletter\ifx\SetFigFont\undefined%
\gdef\SetFigFont#1#2#3#4#5{%
  \reset@font\fontsize{#1}{#2pt}%
  \fontfamily{#3}\fontseries{#4}\fontshape{#5}%
  \selectfont}%
\fi\endgroup%
\begin{picture}(4500,3800)(3500,0)
\hspace{16mm}\includegraphics[width=14cm]{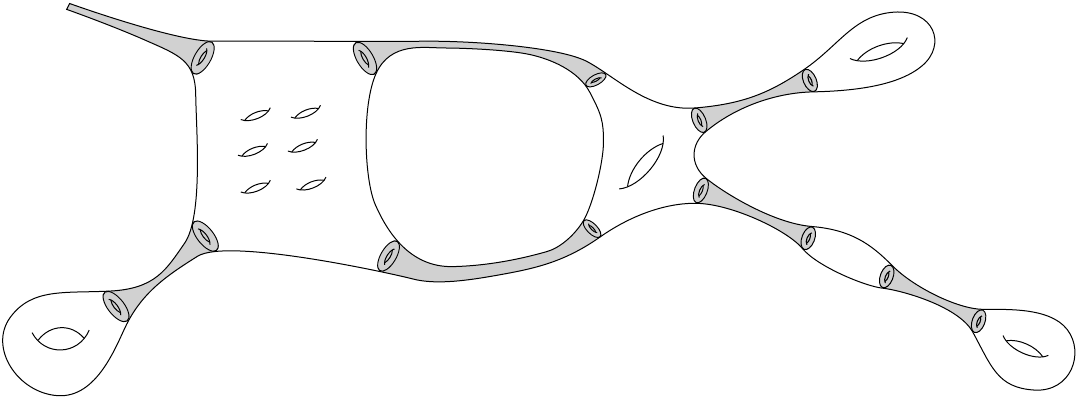}%
\put(-9100,1000){\makebox(0,0)[lb]{\smash{{\SetFigFont{12}{14.4}{\familydefault}{\mddefault}{\updefault}$M_{2,t}$}}}}
\put(-8000,2000){\makebox(0,0)[lb]{\smash{{\SetFigFont{12}{14.4}{\familydefault}{\mddefault}{\updefault}$M_{1,t}$}}}}
\put(-4600,2000){\makebox(0,0)[lb]{\smash{{\SetFigFont{12}{14.4}{\familydefault}{\mddefault}{\updefault}$M_{3,t}$}}}}
\put(-2050,1500){\makebox(0,0)[lb]{\smash{{\SetFigFont{12}{14.4}{\familydefault}{\mddefault}{\updefault}$M_{5,t}$}}}}
\put(-1850,2350){\makebox(0,0)[lb]{\smash{{\SetFigFont{12}{14.4}{\familydefault}{\mddefault}{\updefault}$M_{4,t}$}}}}
\put(-550,870){\makebox(0,0)[lb]{\smash{{\SetFigFont{12}{14.4}{\familydefault}{\mddefault}{\updefault}$M_{6,t}$}}}}
\put(-8200,3300){\makebox(0,0)[lb]{\smash{{\SetFigFont{12}{14.4}{\familydefault}{\mddefault}{\updefault}$E_{1,t}$}}}}
\put(-7800,700){\makebox(0,0)[lb]{\smash{{\SetFigFont{12}{14.4}{\familydefault}{\mddefault}{\updefault}$E_{2,t}$}}}}
\put(-5200,3150){\makebox(0,0)[lb]{\smash{{\SetFigFont{12}{14.4}{\familydefault}{\mddefault}{\updefault}$E_{3,t}$}}}}
\put(-5100,750){\makebox(0,0)[lb]{\smash{{\SetFigFont{12}{14.4}{\familydefault}{\mddefault}{\updefault}$E_{4,t}$}}}}
\put(-3050,1300){\makebox(0,0)[lb]{\smash{{\SetFigFont{12}{14.4}{\familydefault}{\mddefault}{\updefault}$E_{6,t}$}}}}
\put(-3050,2750){\makebox(0,0)[lb]{\smash{{\SetFigFont{12}{14.4}{\familydefault}{\mddefault}{\updefault}$E_{5,t}$}}}}
\put(-1500,550){\makebox(0,0)[lb]{\smash{{\SetFigFont{12}{14.4}{\familydefault}{\mddefault}{\updefault}$E_{7,t}$}}}}
\end{picture}%
\caption{An example for a decomposition of $M_0$.
The subsets $E_{2, t}, \ldots, E_{7,t}$ are diffeomorphic to $T^2 \times I$ and are collapsed along the $T^2$-factor.
The subset $E_{1,t}$ is diffeomorphic to an interval-bundle over the Klein bottle.
The component $M_{1,t}$ is geometrically close to a hyperbolic manifold with cusps and $M_{2,t}, \ldots, M_{6,t}$ collapse along Seifert fibers.
It may happen that the Seifert fibers of $M_{3,t}$ and $M_{4,t}$ are homotopic within $E_{5,t}$ and hence the Seifert fibrations on these two components can be combined to a fibration on $M_{3,t} \cup E_{5,t} \cup M_{4,t}$.
Analogously, $M_{5,t}$ could collapse to an annulus and hence be diffeomorphic to $T^2 \times I$.
In this case, $E_{6,t} \cup M_{5,t} \cup E_{7,t}$ would be diffeomorphic to $T^2 \times I$. 
\label{fig:longtimepicture}}
\end{center}
\end{figure}
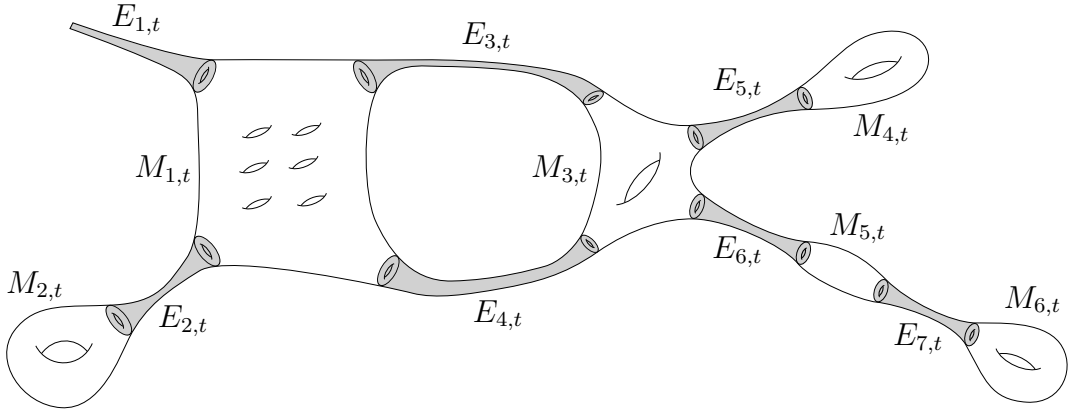
\begin{Theorem} \label{Thm:geombehavior}
Given a Ricci flow with surgery $\MM$ as in Theorem \ref{Thm:MainTheorem-III}, we can find a time $T < \infty$ and a function $\varepsilon : [T, \infty) \to (0, 1)$ with $\lim_{t \to \infty} \varepsilon(t) = 0$ such that the following holds:

The flow $\MM$ has no surgeries past time $T < \infty$, i.e. $\MM$ restricted to the time interval $[T, \infty)$ is a non-singular, long-time existent Ricci flow $(g_t)_{t \in [T, \infty)}$ on some orientable manifold $M$.
Let $M_0 \subset M$ be a component of $M$.
Then $M_0$ is aspherical (i.e. irreducible and not diffeomorphic to a spherical space form) and as $t \to \infty$ the metric $g_t$ on $M_0$ behaves as follows:
\begin{enumerate}[label=(\alph*)]
\item If $M_0$ is diffeomorphic to a hyperbolic manifold, then $\frac14 t^{-1} g_t$ converges to a unique hyperbolic metric on $M_0$.
\item If $M_0$ is diffeomorphic to the $3$-torus $T^3$, then either $g_t$ converges to a flat metric on $M_0$, or $t^{-1/2} \diam_t M_0$ is unbounded on $[T, \infty)$ and for every $t \geq T$ there is a metric $g'_t$ that is $(1+\varepsilon(t))$-bilipschitz close to $g_t$, $t^{-1} g'_t$ is $\varepsilon (t)$-close to $t^{-1} g_t$ in the $C^{[\varepsilon^{-1} (t)]}$-sense\footnote{By this we mean we mean that $\Vert t^{-1} (g_t - g'_t) \Vert_{C^{[\varepsilon^{-1}(t)]}} < \varepsilon (t)$, with respect to the metric $t^{-1} g_t$.} and $g'_t$ is invariant under a free $T^2$-action on $M_0$.
The orbits of this action have diameter $<  \varepsilon(t) \sqrt{t}$ and are the fibers of a $T^2$-fibration over a circle.

If $M_0$ is diffeomorphic to a quotient of $T^3$, then the same statement holds after passing to a finite cover.
\item If $M_0$ is diffeomorphic to a Heisenberg manifold $\Nil^3$, then for every $t \geq T$ there is a metric $g'_t$ such that $g_t$ is $(1+\varepsilon(t))$-bilipschitz close to $g_t$, $t^{-1} g'_t$ is $\varepsilon (t)$-close to $t^{-1} g_t$ in the $C^{[\varepsilon^{-1} (t)]}$-sense and such that the following holds:

Either we have $\diam_t M_0  < \varepsilon(t) \sqrt{t}$ for all $t \geq T$ and the metrics $g'_t$ are isometric to quotients of left-invariant metrics on $\Nil^3$, or the following holds:
The normalized diameter $t^{-1/2} \diam_t M_0$ is unbounded on $[T, \infty)$ and for every $t \geq T$ we can express $M_0$ as a $T^2$-bundle over a circle such that in a fibered neighborhood of every $T^2$-fiber there is a free $T^2$-action that is isometric with respect to $g'_t$ and whose orbits are the $T^2$-fibers.
Moreover, the diameter of each $T^2$-fiber with respect to $g_t$ is  $<  \varepsilon(t) \sqrt{t}$.

If $M_0$ is diffeomorphic to a quotient of $\Nil^3$, then the same statement holds after passing to a finite cover.
\item If $M_0$ is diffeomorphic to a Solv manifold, then for every $t \geq T$ there is a metric $g'_t$ such that $g_t$ is $(1+\varepsilon(t))$-bilipschitz close to $g_t$, $t^{-1} g'_t$ is $\varepsilon (t)$-close to $t^{-1} g_t$ in the $C^{[\varepsilon^{-1} (t)]}$-sense and such that the following holds:
For every $t \geq T$ we can express $M_0$ as a $T^2$-bundle over a circle such that in a fibered neighborhood of every $T^2$-fiber there is a free $T^2$-action that is isometric with respect to $g'_t$ and whose orbits are the $T^2$-fibers.
Moreover, the diameter of each $T^2$-fiber with respect to $g_t$ is  $<  \varepsilon(t) \sqrt{t}$.

If $t^{-1/2} \diam_t M_0$ stays bounded on $[T, \infty)$, then for all $t \geq T$ the metric $g'_t$ is isometric to a quotient of a left-invariant metric on the Solv Lie-group.

If $M_0$ is diffeomorphic to a quotient of the Solv manifold, then the same statement holds after passing to a finite cover.
\item In all other cases we have the following picture:
There is a constant $A_0 < \infty$, which only depends on the topology of $M_0$, and for every $\mu > 0$, there are constants\footnote{Note that, unlike $A_0$, the constants $a, B$ and $T_0$ may depend not only on the topology of $M_0$, but on the geometry of $\MM$.} $a(\mu) > 0$ and $B(\mu), T_0(\mu) < \infty$ such that:

Let $(H_1, g_{\textnormal{hyp},1}), \ldots, (H_p, g_{\textnormal{hyp},p})$ be the hyperbolic manifolds (of finite volume) whose underlying topological manifolds occur as hyperbolic pieces in the geometric decomposition of $M_0$.
For each $j = 1, \ldots, p$ and sufficiently small $b > 0$ denote by $H^{(b)}_j$ the manifold that arises from $H_j$ by chopping off the cusps along horospherical, cross-sectional $2$-tori of area $b$.

Then for each $t \geq T_0 (\mu)$ we can find a metric $g'_t$ on $M_0$ that is $(1+\mu)$-bilipschitz close to $g_t$ and $\mu$-close to $g_t$ in the $C^{[\mu^{-1}]}$-sense.
Moreover, there are finitely many, pairwise disjoint subsets $E_{1,t}, \ldots, E_{m_t,t} \subset M_0$ such that the following holds:
Let $M_{1,t}, \ldots, \linebreak[1] M_{k_t, t} \subset M_0$ be the closures of components $M_0 \setminus (E_{1,t} \cup \ldots \cup E_{m_t,t})$.
Then
\begin{enumerate}[label=(e\arabic*)]
\item Each $E_{i,t}$ is diffeomorphic to $I \times T^2$ or to a twisted interval bundle over the Klein bottle, $\Klein^2 \td\times I$.
The (generic) $T^2$-fibers of $E_{i,t}$ are incompressible (i.e. $\pi_1$-injective) in $M_0$.
\item We have $k_t \geq p$ and after possibly relabeling the $M_{j,t}$ we have:
For each $j = 1, \ldots, p$ the interior of $M_{j,t}$ is diffeomorphic to $H_j$.
For each $j = p+1, \ldots, k_t$ the component $M_{j,t}$ admits a Seifert fibration $p_{j,t} : M_{j,t} \to \Sigma_{j,t}$, where $\Sigma_{j,t}$ is an orbifold with boundary whose singularities are of cone type.
\item For all $i = 1, \ldots, m_t$ we have $\diam_t E_{i,t} > \mu^{-1} \sqrt{t}$.
\item For all $j = 1, \ldots, k_t$ we have $\diam_t M_{j,t} < B(\mu) \sqrt{t}$.
\item For each $E_{i,t}$ that is diffeomorphic to $T^2 \times I$ there is a diffeomorphism $\Phi_{E_{i,t}} : T^2 \times I \to E_{i,t}$ such that $\Phi^*_{E_{i,t}} g'_t$ is invariant under the $T^2$-action on the first factor. 
Moreover, the orbits of this action have diameter $< \mu \sqrt{t}$ and second fundamental form $< B(\mu) t^{-1/2}$ with respect to $\Phi^*_{E_{i,t}} g'_t$ and $\Phi^*_{E_{i,t}} g_t$.

If $E_{i,t} \approx \Klein \td\times I$, then the same statements holds for the double cover that is diffeomorphic to $T^2 \times I$.
\item For each $j = 1, \ldots, p$ there is a diffeomorphism $\Phi_{j,t} : H_{j}^{(\varepsilon (t))} \to M_{j,t}$ such that $\frac14 t^{-1} \Phi^*_{j,t} g_t$ is $\varepsilon(t)$-close to the hyperbolic metric $g_{\textnormal{hyp},j}$ on $H_{j}^{(\varepsilon(t))}$ in the $C^{[\varepsilon^{-1}(t)]}$-sense.

\item For each $j = p+1, \ldots, k_t$, the fibers of the Seifert fibration on $M_{j,t}$ have diameter $< \varepsilon (t) \sqrt{t}$ and are orbits of an $S^1$-action on $M_{j,t}$ that is isometric with respect to $g'_t$.
There is an orbifold metric $g''_{j,t}$ on $\Sigma_{j,t}$ such that the projection map $p_{j,t} : (M_{j,t}, g'_t) \to (\Sigma_{j,t}, g''_t)$ is a submersion.
The geodesic curvature of the Seifert fibers on $M_{j,t}$ is bounded by $B(\mu) t^{-1/2}$ with respect to $g'_t$ and $g_t$ and the curvature on $(\Sigma_{j,t}, g''_t)$ is bounded by $B(\mu) t^{-1}$.
Moreover, on $M_{j,t}$ the metrics $g_t$ and $g'_t$ are even $\varepsilon(t)$-close in the $C^{[\mu^{-1}]}$-sense.

\item We have the area bounds
\[ \qquad \qquad \area (\Sigma_{j,t}, g''_t) > a(\mu) t \qquad \text{for all} \qquad j = p+1, \ldots, k_t \]
and
\[ \area (\Sigma_{p+1,t}, g''_t) + \ldots + \area (\Sigma_{k_t, t}, g''_t) < A_0 t. \]
\item Every component of $M_0 \setminus \Int (M_{1, t} \cup \ldots \cup M_{p,t})$ that is diffeomorphic to $T^2 \times I$ is equal to some $E_{i,t}$.
\end{enumerate}
\end{enumerate}
\end{Theorem} 

The subsets $E_{1,t}, \ldots, E_{m_t, t}$ that are diffeomorphic to $T^2 \times I$ can be interpreted as tubular neighborhoods of the incompressible $2$-tori of a geometric decomposition of $M_0$.
This geometric decomposition is not necessarily minimal (compare again with Definition \ref{Def:geomdec} of the third paper in this series, \cite{Bamler-LT-topology}).
For example, we did not exclude the possibility that there is a component $E_{i,t}$ that is diffeomorphic to $T^2 \times I$ and that has the property that the Seifert fibers coming from the two adjacent $M_{j,t}$ are homotopic to each other within $E_{i,t}$.
Such a component $E_{i,t}$ would correspond to a redundant torus in the geometric decomposition, because the Seifert fibrations on the two adjacent $M_{j,t}$ could be extended (topologically) to a Seifert fibration on the union of those two $M_{j,t}$ and the connecting $E_{i,t}$.
Another possibility that we did not exclude in Theorem \ref{Thm:geombehavior} is that one of the $M_{j,t}$ is a (non-singular) Seifert fibration over an annulus and therefore diffeomorphic to $T^2 \times I$.
In this case, the $T^2$-fibration of the two adjacent $E_{i,t}$ can be extended topologically onto $M_{j,t}$ in an analogous way, making the $2$-torus corresponding to one of those $E_{i,t}$ redundant.
In both examples, the extension process would simplify the decomposition of $M_0$.
This simplification, however, is then not reflected by the metric $g_t$.
The previous theorem makes no statement about whether such pieces $E_{i,t}$ that correspond to redundant tori in the geometric decomposition can occur.
See again Figure \ref{fig:longtimepicture} for an illustration of these two possibilities.

For a more concrete example consider the case in which $M_0 \approx \Sigma \times S^1$, where $\Sigma$ is a surface of genus $g \geq 2$.
Then a priori, there is no bound on $m_t$ that only depends on the topology of $M_0$, since a geometric decomposition can be induced from cutting $\Sigma$ along arbitrarily many pairwise disjoint, embedded, non-contractible loops.
If the number of these loops is large enough, then some components in their complement are annuli and hence some of the $M_{j,t}$ are diffeomorphic to $T^2 \times I$.
So the following question arises naturally:

\begin{Question}
In part (e) of Theorem \ref{Thm:geombehavior}, can we choose $E_{1,t}, \ldots, E_{m_t, t}$ such that none of the components $M_{j,t}$ are diffeomorphic to $T^2 \times I$ or $\Klein \td\times I$?
\end{Question}

If that was the case, then $m_t$ would be uniformly bounded in terms of the topology of $M_0$.
More generally, we may ask:

\begin{Question}
Can $E_{1,t}, \ldots, E_{m_t, t}$ be chosen such that the corresponding geometric decomposition is minimal?
\end{Question}

Note that minimal geometric decompositions are unique up to isotopy (cf \cite[sec \ref{sec:3dtopology}]{Bamler-LT-topology}).
An affirmative answer to this question would imply that the Seifert fibers on either side of each $E_{j,t}$ are not homotopic to each other.
Moreover, $m_t$ would be a constant depending only on the topology of $M_0$.
In the setting in which $M_0$ consists of a single geometric component that is not flat, nil or solv, the previous question is equivalent to the fact that $m_t = 0$ and to the following question (via Theorem \ref{Thm:geombehavior}(e3), (e4)):

\begin{Question}
Assume that the geometric decomposition of $M_0$ consists of a single component.
Is there a constant $C<\infty$ such that $\diam_t M_0 < C \sqrt{t}$ for all $t \geq T$?
\end{Question}

If this diameter bound was known, then by the work of Lott (\cite{LottDimRed}) we can understand the subsequential Gromov-Hausdorff limits of $(M_0, t^{-1} g_t)$ as $t \to \infty$ as well as the Gromov-Hausdorff limits of the universal covers of $(M_0, t^{-1} g_t)$.
More generally, we may ask:

\begin{Question}
Assume that we are in case (e) of Theorem \ref{Thm:geombehavior}.
Do the submersion metrics $t^{-1} g''_t$ on each base orbifold $\Sigma_{j,t}$ limit to a certain standard metric?
Or more generally: do the metrics $t^{-1} g_t$ on larger and larger neighborhoods around $M_{j,t}$ collapse to certain standard geometries?
When are these geometries unique, i.e. when do they only depend on the topology of $M_0$.
\end{Question}

Natural candidates for such  standard geometries would be the hyperbolic surfaces of curvature $-\frac12$.
Motivated by this picture, a further question would be the following:

\begin{Question}
Is there a function $\varepsilon : [T, \infty) \to (0,\infty)$ such that $\lim_{t \to \infty} \varepsilon(t) = 0$ and such that for all $t \geq T$
\[ - \big( \tfrac12 + \varepsilon(t) \big) t^{-1} < \sec_t < \big( \tfrac1{4} +  \varepsilon(t) \big) t^{-1}  \quad \text{and} \quad - \big( \tfrac12 + \varepsilon(t) \big) t^{-1}  < \Ric_t < \big( \tfrac16 + \varepsilon(t) \big) t^{-1}  \]
on $M_0$?
\end{Question}
(The two lower bounds are realized by the geometric models $\IH^2 \times \IR$, $PSL (2, \IR)$ and $Sol$, the upper sectional curvature bound is realized by $Sol$ and the upper Ricci curvature bound is realized by $Nil$.)

Finally, we may still ask:

\begin{Question}
Does the metric $t^{-1} g_t$ converge pointwise to a possibly singular metric $g_\infty$ as $t \to \infty$.
\end{Question}

Note that all these questions reduce to questions about \emph{non-singular} Ricci flows $(g_t)_{t \in [0,\infty)}$ in the wake of Theorem \ref{Thm:MainTheorem-III}.

\subsection{Outline of the proofs} \label{subsec:outlineofintroduction}
We will now give a brief outline of the proofs of Theorems \ref{Thm:geombehavior} and \ref{Thm:MainTheorem-III}.
More detailed explanations of specific aspects of the proofs can be found in the introductions of the subsequent four papers.

The most important finding of these four papers is the curvature bound $|{\Rm_t}| < Ct^{-1}$ for large $t$.
Using this bound, it is possible to rule out the existence of surgeries at large times, since surgeries only occur where the curvature blows up.
Moreover, the geometric characterization of Theorem \ref{Thm:geombehavior} follows from looking at the proof of this curvature bound.
It will turn out that in order to establish said curvature bound, the existence of surgeries does not create any major issues, apart from several technical difficulties.
So in this outline we will restrict ourselves to the case in which the given Ricci flow is non-singular, i.e. it is given by a smooth family of metrics $(g_t)_{t \in [0, \infty)}$ on $M$.
In other words, in the following we will sketch the proof of Corollary \ref{Cor:nonsingflow}.

For this, we have to recall an important result of Perelman, which in this outline we will refer to as the ``Key Lemma''.
Define for every point $(x,t) \in M \times [0, \infty)$ in space-time the scale $\rho(x,t) > 0$ as follows
\[ \rho (x,t) := \sup \big\{ r > 0 \;\; : \;\; \sec_t \geq - r^{-2} \; \text{on} \; B(x,t, r) \big\}. \]
The condition after the colon reads that the sectional curvature on a ball at time-$t$ radius $r$ around $x$ is bounded at time $t$ from below by $-r^{-2}$.
Equivalently, we could say that $\rho(x,t)$ is the maximal radius such that if we rescale the ball $B(x,t, \rho(x,t))$ to have radius $1$, then the sectional curvature on this ball is bounded from below by $-1$.
If the sectional curvature on the component of $M$ in which $x$ lies is non-negative at time $t$, then we write $\rho(x,t) = \infty$.
Now Perelman's Key Lemma can be phrased as follows:

\begin{KeyLemma}
For every $w > 0$ there are constants $\ov{\rho}(w) > 0$ and $K(w) < \infty$ such that if
\begin{equation} \label{eq:rhoballisnotcollapsed} \vol_t B(x,t, \rho(x,t)) > w \rho^3(x,t), \end{equation}
then
\[ \rho(x,t) > \ov\rho (w) \sqrt{t} \qquad \text{and} \qquad |{\Rm}|(x,t) < \frac{K(w)}{t}. \]
\end{KeyLemma}

It is even possible to obtain the same curvature bound on the ball $B(x,t, \linebreak[1] \ov\rho (w) \sqrt{t})$ and for times of the time-interval $[(1-\tau(w))t, t]$ for some uniform $\tau(w) > 0$.
Motivated by this Key Lemma, Perelman decomposes the manifold $M$ into a \emph{thick} part $M_{\textnormal{thick}}(t)$ and a \emph{thin} part $M_{\textnormal{thin}}(t)$ for every time $t$:
\[ M = M_{\textnormal{thick}} (t) \; \dotcup \; M_{\textnormal{thin}} (t) \]
The thick part roughly consists of all points $x \in M$ that at time $t$ satisfy (\ref{eq:rhoballisnotcollapsed}) for a suitable $w$ and the thin part is the complement of the thick part.
So on $M_{\textnormal{thick}}(t)$ the curvature is bounded by $K(w) t^{-1}$ and $M_{\textnormal{thin}}(t)$ is locally collapsed at the scale $\rho$.
Using the curvature bound on the thick part, Perelman could show that for sufficiently large $t$, the metric on $M_{\textnormal{thick}}(t)$ is close to a hyperbolic metric.
On the other hand, using collapsing theory (\cite{ShioyaYamaguchi}, \cite{MorganTian}, \cite{KLcollapse}) it is possible to decompose $M_{\textnormal{thin}}(t)$ into regions that exhibit different collapsing behaviors at the scale $\rho$.
Generically, those collapses can occur along $S^1$, $T^2$ or $S^2$ fibers or $M$ could be globally collapsed.
The decomposition of $M_{\textnormal{thin}}(t)$ into regions of different collapsing behaviors arises from cutting $M_{\textnormal{thin}} (t)$ along embedded $2$-spheres or $2$-tori.
Note that some of these $2$-spheres or $2$-tori could be contractible or compressible in $M$.
So, a priori, this decomposition can be very different from, and far more complex than the geometric decomposition of $M$ (see \cite[Definition \ref{Def:geomdec}]{Bamler-LT-topology}).
Further topological arguments are needed in order to reorganize this decomposition into a geometric decomposition of $M$.

For the Main Theorem \ref{Thm:MainTheorem-III}, it suffices to establish the desired curvature bound on the thin part of $M$ since this is the part of the manifold where Perelman's Key Lemma fails.
In order to achieve this bound, we make use of the observation that Perelman's Key Lemma continues to hold if we pass to the universal cover.
By this we mean the following: Consider the universal cover $\pi : \td{M} \to M$ of $M$ and pull back the family of metrics $g_t$ to $\td{M}$ via the projection map $\pi$.
Then these pull-backs $\pi^* g_t$ still satisfy the Ricci flow equation.
It turns out that Perelman's proof also works in this (possibly non-compact) setting.
Now consider for every point $x \in M$ one of its lifts $\td{x} \in \td{M}$ (i.e. $\pi (\td{x}) = x$) and look at the ball $B^{\td{M}} ( \td{x}, t, \rho(x,t) )$ around $\td{x}$ in $(\td{M}, \pi^* g_t)$.
The volume of this ball is not smaller than the volume of the ball $B(x,t, \rho(x,t))$ in $M$.
Hence we can generalize Perelman's Key Lemma as follows:

\begin{GeneralizedKeyLemma}
For every $w > 0$ there are constants $\ov{\rho}(w) > 0$ and $K(w) < \infty$ such that if
\begin{equation} \label{eq:rhoballnotcollapsedinuniversalcover} \vol_t B^{\td{M}} ( \td{x}, t, \rho(x,t) ) > w \rho^3(x,t), \end{equation}
then
\[ \rho(x,t) > \ov\rho (w) \sqrt{t} \qquad \text{and} \qquad |{\Rm}|(x,t) < \frac{K(w)}{t}. \]
\end{GeneralizedKeyLemma}

In subsection \ref{subsec:geometricconsequences} of \cite{Bamler-LT-main} will see that the bound (\ref{eq:rhoballnotcollapsedinuniversalcover}) can be guaranteed for a suitable $w$, whenever the metric around $x$ is either non-collapsed at scale $\rho(x,t)$ or whenever it is collapsed at scale $\rho(x,t)$ along \emph{incompressible} $S^1$ or $T^2$-fibers.
Recall that by ``incompressible'' we mean that the fundamental group of the fibers injects into the fundamental group of $M$.
From now on we will call regions of $M_{\textnormal{thin}} (t)$ where such a collapse occurs \emph{good} and the remaining regions \emph{bad}.
So we obtain a decomposition
\[ M_{\textnormal{thin}}(t) = M_{\textnormal{good}}(t) \; \dotcup \; M_{\textnormal{bad}}(t). \]
Summarizing our results, we can say that we have established the curvature bound $|{\Rm_t}| < C t^{-1}$ on $M_{\textnormal{thick}}(t) \cup M_{\textnormal{good}}(t) = M \setminus M_{\textnormal{bad}} (t)$.

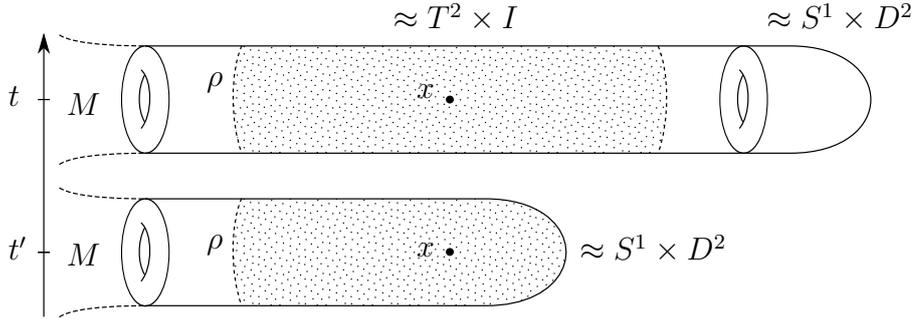
\begin{figure}[t] 
\begin{center}
\setlength{\unitlength}{2863sp}%
\begingroup\makeatletter\ifx\SetFigFont\undefined%
\gdef\SetFigFont#1#2#3#4#5{%
  \reset@font\fontsize{#1}{#2pt}%
  \fontfamily{#3}\fontseries{#4}\fontshape{#5}%
  \selectfont}%
\fi\endgroup%
\begin{picture}(2500,3000)(3500,0)
\hspace{16mm}\includegraphics[width=11cm]{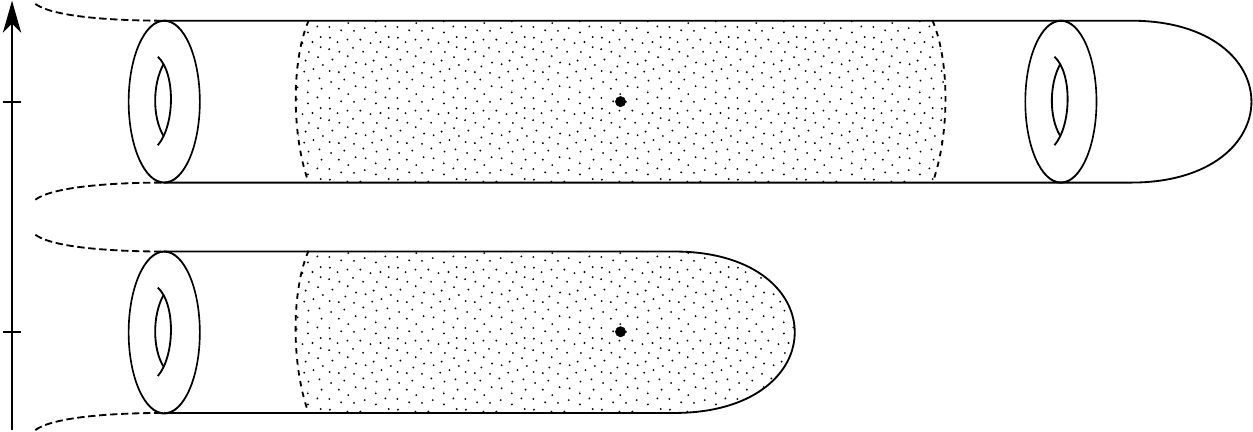}%
\put(-3900,530){\makebox(0,0)[lb]{\smash{{\SetFigFont{12}{14.4}{\familydefault}{\mddefault}{\updefault}$x$}}}}
\put(-6900,1750){\makebox(0,0)[lb]{\smash{{\SetFigFont{12}{14.4}{\familydefault}{\mddefault}{\updefault}$M$}}}}
\put(-6900,450){\makebox(0,0)[lb]{\smash{{\SetFigFont{12}{14.4}{\familydefault}{\mddefault}{\updefault}$M$}}}}
\put(-3900,1900){\makebox(0,0)[lb]{\smash{{\SetFigFont{12}{14.4}{\familydefault}{\mddefault}{\updefault}$x$}}}}
\put(-5700,2000){\makebox(0,0)[lb]{\smash{{\SetFigFont{12}{14.4}{\familydefault}{\mddefault}{\updefault}$\rho$}}}}
\put(-5700,600){\makebox(0,0)[lb]{\smash{{\SetFigFont{12}{14.4}{\familydefault}{\mddefault}{\updefault}$\rho$}}}}
\put(-2500,500){\makebox(0,0)[lb]{\smash{{\SetFigFont{12}{14.4}{\familydefault}{\mddefault}{\updefault}$\approx S^1 \times D^2$}}}}
\put(-4100,2500){\makebox(0,0)[lb]{\smash{{\SetFigFont{12}{14.4}{\familydefault}{\mddefault}{\updefault}$\approx T^2 \times I$}}}}
\put(-900,2500){\makebox(0,0)[lb]{\smash{{\SetFigFont{12}{14.4}{\familydefault}{\mddefault}{\updefault}$\approx S^1 \times D^2$}}}}
\put(-7400,1830){\makebox(0,0)[lb]{\smash{{\SetFigFont{12}{14.4}{\familydefault}{\mddefault}{\updefault}$t$}}}}
\put(-7400,500){\makebox(0,0)[lb]{\smash{{\SetFigFont{12}{14.4}{\familydefault}{\mddefault}{\updefault}$t'$}}}}
\end{picture}%
\caption{In this example the universal cover of the ball $B(x,t,\rho)$, which is contained in the region that is diffeomorphic to $T^2 \times I$, is larger than the ball around a lift $\td{x}$ in the universal cover of $\td{M}$, because some loops that are non-contractible within $B(x,t, \rho)$ are contractible in the region that is diffeomorphic to $S^1 \times D^2$.
Perelman's estimates cannot be localized easily in this case, since $B(x,t', \rho)$ may include this region at some time $t' < t$.\label{fig:expandingsolidtorus}}
\end{center}
\end{figure}
Note, that the ball $B^{\td{M}} ( \td{x}, t, \rho) $ in the universal cover of $M$ is in general not equal to the universal cover of the ball $B(x, t, \rho)$, which we would denote by $\td{B(x,t, \rho)}$.
The volume of a $\rho$-ball in this cover is in general even bigger than the volume of $B^{\td{M}} ( \td{x}, t, \rho) )$.
But unfortunately, Perelman's Key Lemma does in general not seem to hold if we replace the volume in (\ref{eq:rhoballnotcollapsedinuniversalcover}) by the volume of a $\rho$-ball in the universal cover of $B(x,t,\rho)$.
This is why we need the collapse in $M_{\textnormal{good}}(t)$ to occur along fibers that are \emph{incompressible in $M$}.
We explain briefly, why it seems unlikely to prove such a curvature bound assuming this more general volume bound:
In order to prove such a bound, we would have to pass to a (local) cover of $B(x,t,\rho)$.
The proof of Perelman's Key Lemma can only be carried out in such a local cover, if we can ensure that all its estimates only take place very close to $x$.
However, these estimates involve the metric at slightly earlier times and due to the lack of a priori curvature bounds, we have no control on the distance distortion under the Ricci flow.
So points that are very close to $x$ at some time $t' < t$ can lie outside $B(x,t,\rho)$ at time $t$.
Figure \ref{fig:expandingsolidtorus} illustrates this problem.
Here the ball $B(x,t,\rho)$ has fundamental group $\IZ^2$, as does the region that is diffeomorphic to $T^2 \times I$, but this region is contained in a union of two regions, which as a whole is diffeomorphic to a solid torus $S^1 \times D^2$ and hence has fundamental group $\IZ$.
In other words, the homotopy classes corresponding to one of the $\IZ$ factors of the fundamental group of the left region are ``destroyed'' by the right region.
Since we don't know whether the right region is disjoint from the ball $B(x,t',\rho)$ at the earlier time $t'$, we cannot pass to an appropriate local cover.
Thus a further generalization of Perelman's Key Lemma in this broad setting seems unlikely.
It will become important later on, however, that if we can exclude such a behavior, then it is in fact possible to localize Perelman's arguments and prove an even more general Key Lemma in certain settings.

Let us now return to our analysis of $M_{\textnormal{thin}} (t)$.
Recall that in order to prove the Main Theorem, it remains to establish the desired curvature bound on $M_{\textnormal{bad}}(t)$.
In the following paragraphs we will roughly sketch how this bound is obtained.
For more details we refer to the introductions of the following papers, particularly of the last paper of this series, \cite{Bamler-LT-main}.

\begin{figure}[t] 
\begin{center}
\setlength{\unitlength}{2863sp}%
\begingroup\makeatletter\ifx\SetFigFont\undefined%
\gdef\SetFigFont#1#2#3#4#5{%
  \reset@font\fontsize{#1}{#2pt}%
  \fontfamily{#3}\fontseries{#4}\fontshape{#5}%
  \selectfont}%
\fi\endgroup%
\begin{picture}(4000,5300)(3500,0)
\hspace{16mm}\includegraphics[width=14cm]{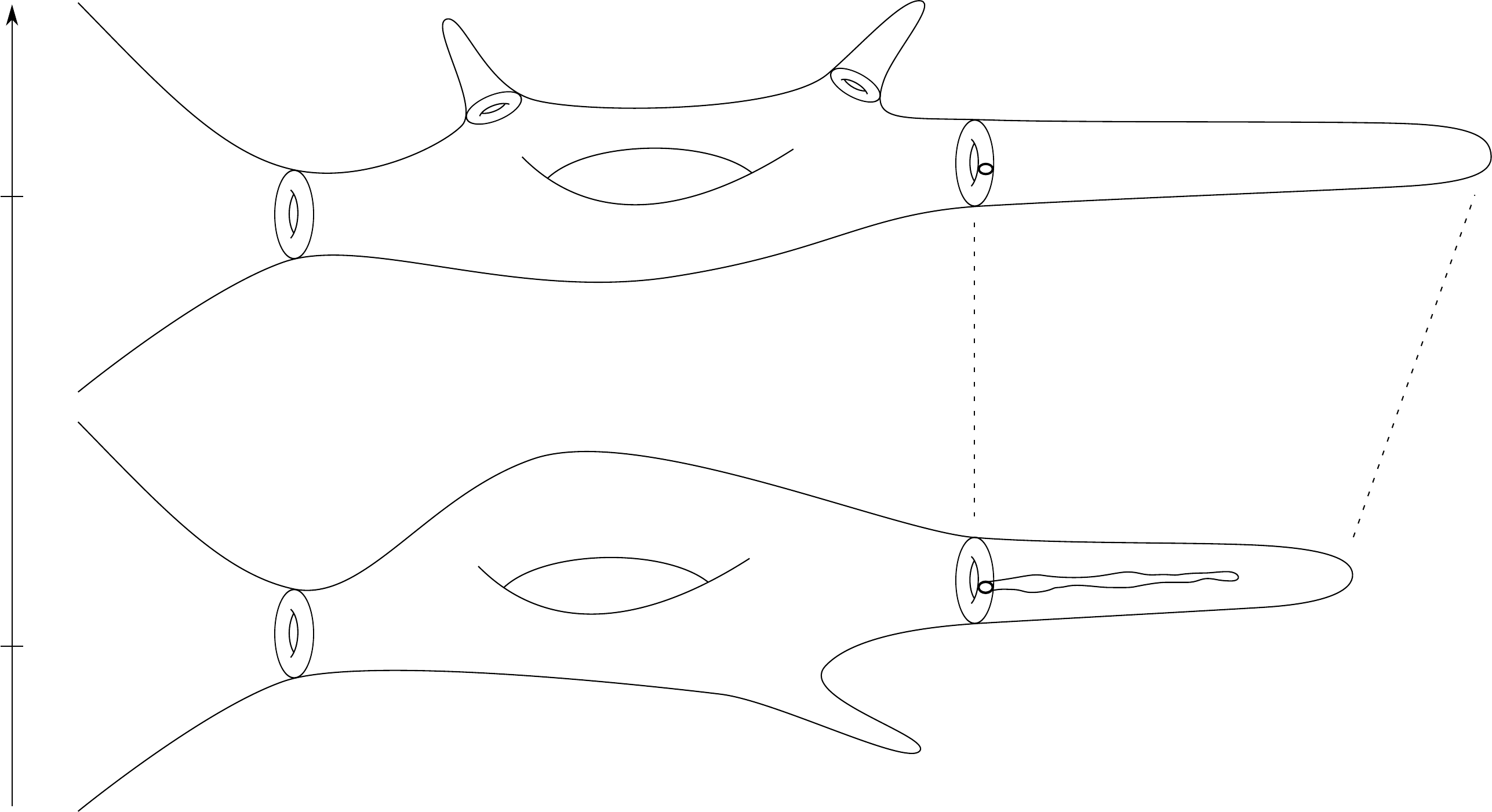}%
\put(-6830,4600){\makebox(0,0)[lb]{\smash{{\SetFigFont{12}{14.4}{\familydefault}{\mddefault}{\updefault}$S_{1,t}$}}}}
\put(-4280,4800){\makebox(0,0)[lb]{\smash{{\SetFigFont{12}{14.4}{\familydefault}{\mddefault}{\updefault}$S_{2,t}$}}}}
\put(-6700,3500){\makebox(0,0)[lb]{\smash{{\SetFigFont{12}{14.4}{\familydefault}{\mddefault}{\updefault}$M_{\textnormal{thin}} (t)$}}}}
\put(-8700,3500){\makebox(0,0)[lb]{\smash{{\SetFigFont{12}{14.4}{\familydefault}{\mddefault}{\updefault}$M_{\textnormal{thick}} (t)$}}}}
\put(-6700,1000){\makebox(0,0)[lb]{\smash{{\SetFigFont{12}{14.4}{\familydefault}{\mddefault}{\updefault}$M_{\textnormal{thin}} (t)$}}}}
\put(-8700,1000){\makebox(0,0)[lb]{\smash{{\SetFigFont{12}{14.4}{\familydefault}{\mddefault}{\updefault}$M_{\textnormal{thick}} (t)$}}}}
\put(-2950,3900){\makebox(0,0)[lb]{\smash{{\SetFigFont{12}{14.4}{\familydefault}{\mddefault}{\updefault}$\gamma_{3,t}$}}}}
\put(-3020,1520){\makebox(0,0)[lb]{\smash{{\SetFigFont{12}{14.4}{\familydefault}{\mddefault}{\updefault}$\gamma_{3,t}$}}}}
\put(-1485,1390){\makebox(0,0)[lb]{\smash{{\SetFigFont{12}{14.4}{\familydefault}{\mddefault}{\updefault}$h'_{3,t}$}}}}
\put(-1600,3900){\makebox(0,0)[lb]{\smash{{\SetFigFont{12}{14.4}{\familydefault}{\mddefault}{\updefault}$S_{3,t}$}}}}
\put(-9300,3700){\makebox(0,0)[lb]{\smash{{\SetFigFont{12}{14.4}{\familydefault}{\mddefault}{\updefault}$t$}}}}
\put(-9550,950){\makebox(0,0)[lb]{\smash{{\SetFigFont{12}{14.4}{\familydefault}{\mddefault}{\updefault}$\alpha_0 t$}}}}
\end{picture}%
\caption{The curvature is bounded outside the solid tori $S_{1,t}, S_{2,t}, S_{3,t}$.
The solid tori $S_{1,t}, S_{2,t}$ have time-$t$ diameter $< D_0 \sqrt{t}$ and the curvature is bounded on $S_{1,t} \cup S_{2,t}$ by $K_0 t^{-1}$.
The diameter of $S_{3,t}$ is $> D_0 \sqrt{t}$ and $S_{3,t}$ persists until time $\alpha_0 t < t$.
We can find a loop $\gamma_{3,t} \subset \partial S_{3,t}$ that is short and geometrically controlled on the time-interval $[\alpha_0 t, t]$.
At time $\alpha_0 t$ this loop spans a disk $h'_{3,t} : D^2 \to M$ of time-$(\alpha_0 t)$ area $< A_0 \alpha_0 t$.\label{fig:picturewithSi}}
\end{center}
\end{figure}
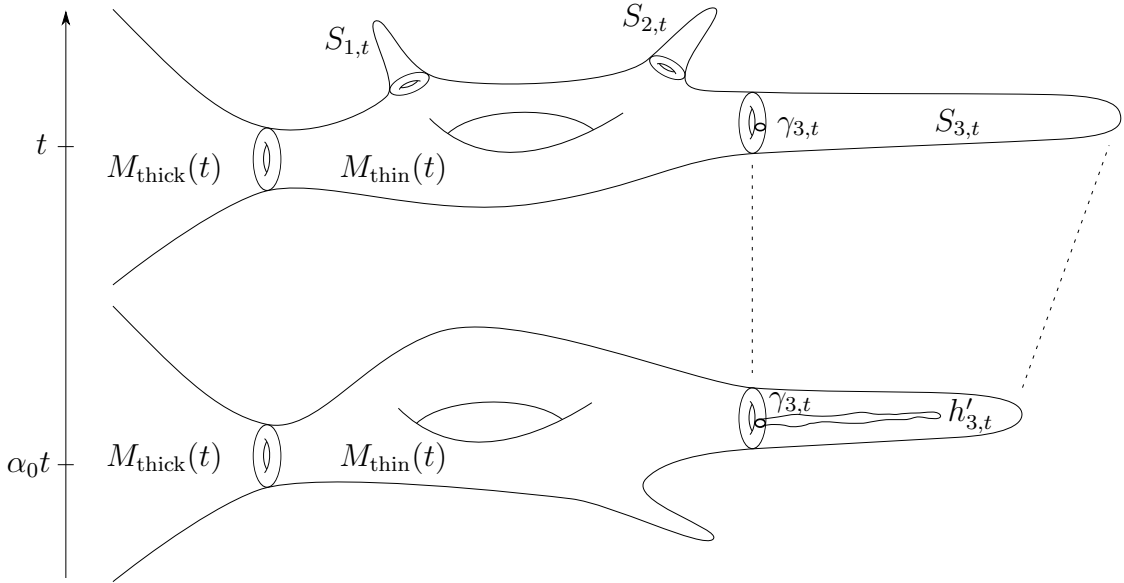
First, we analyze the topology of the decomposition of $M_{\textnormal{thin}}(t)$ to understand its fragmentation into good and bad parts.
We will learn that $M_{\textnormal{bad}}(t)$ can be covered by pairwise disjoint regions that are either diffeomorphic to $T^2 \times I$ or to solid tori $S^1 \times D^2$.
On those regions that are diffeomorphic to $T^2 \times I$ it is in fact possible to localize Perelman's Key Lemma.
The reason why we can do this comes from the fact that the part of $M$ that ``destroys'' certain homotopy classes of the $T^2 \times I$ regions stays sufficiently far away from this region for a short but uniform time-interval.
So we conclude that there are pairwise disjoint, embedded solid tori $S_{1,t}, \ldots, S_{m_t, t} \subset M$, $S_{i,t} \approx S^1 \times D^2$ such that
\[ |{\Rm_t}| < K t^{-1} \qquad \text{on} \qquad M \setminus (S_{1,t} \cup \ldots \cup S_{m_t, t}). \]
Note that the number and position of the solid tori $S_{1,t}, \ldots, S_{m_t, t}$ depends on $t$.

It now remains to analyze the solid tori $S_{i,t}$.
This analysis requires further generalizations of Perelman's arguments in the collapsing case, e.g. bounded curvature at bounded distance results.
A major goal of this analysis is to show that a behavior as illustrated in Figure \ref{fig:expandingsolidtorus} cannot occur.
This will allow us to apply a localized version of Perelman's Key Lemma yielding curvature control on large collar neighborhoods of those $S_{i,t}$ whose diameter is large.
We refer to the introduction of the last paper of this series, \cite{Bamler-LT-main}, for further details on the analysis of the solid tori $S_{i,t}$.

A major ingredient that is used in this analysis is the fact that, in most cases, the solid tori $S_{i,t}$ admit compressing disks of bounded area.
By this we mean that for every solid torus $S_{i,t}$ we can find a map $h_{i,t} : D^2 \to M$ with the following properties: $h_{i,t} (\partial D^2) \subset \partial S_{i,t}$ and $h_{i,t} |_{\partial D^2}$ parameterizes a loop that is non-contractible in $\partial S_{i,t}$ (but contractible in $M$).
Moreover, the time-$t$ area of $h_{i,t}$, is $< A_0 t$ for some uniform constant $A_0 < \infty$ that only depends on the topology of $M$.
These compressing disks roughly arise from the intersection of $S_{i,t}$ with certain ``minimal simplicial complexes'', whose area is bounded by $A t$ for some other uniform $A < \infty$.
Note that both the area bound for these ``minimal simplicial complexes'' and the extraction of said compressing disks are two non-trivial steps in our proof and occupy most of the second and third paper, \cite{Bamler-LT-simpcx} and \cite{Bamler-LT-topology}.
Furthermore, unfortunately, the extraction of the compressing disk seems to fail in the special case in which $M$ is topologically a quotient of a $2$-torus bundle over a circle.
In this case, we employ a different argument, which makes use of the special topology of $M$.

Taking all these facts into account, the analysis of the $S_{i,t}$ yields the following conclusion for each $i$:
Either the diameter of $S_{i,t}$ is bounded by a constant of the form $D_0 \sqrt{t}$ and the curvature on $S_{i,t}$ is bounded by a constant of the form $K_0 t^{-1}$ or the diameter of $S_{i,t}$ is larger than $D_0 \sqrt{t}$ and the solid torus $S_{i,t}$ persists in a certain sense if we go backwards in the flow on a long time-interval of the form $[\alpha_0 t, t]$.
Moreover, after passing to a smaller solid torus, we find a loop $\gamma_{i,t} \subset \partial S_{i,t}$ that is non-contractible in $\partial S_{i,t}$, but contractible in $M$, that is short and whose geodesic curvature is sufficiently controlled on the whole time-interval $[\alpha_0 t, t]$.
This loop can be chosen such that it spans a disk $h'_{i,t} : D^2 \to M$ of area $< A_1 \alpha_0 t$ at time $\alpha_0 t$, where $A_1$ can be determined from $A_0$ and the topology of $M$.
See Figure \ref{fig:picturewithSi} for an illustration of these two cases.

We can finally rule out the second case using a minimal disk argument, which is due to Hamilton.
This argument implies that the loop $\gamma_{i,t}$ can only be short and have bounded geodesic curvature on a maximal time-interval of the form $[\alpha_0 t, B(A_1) \alpha_0 t]$, where $B(A_1) < \infty$ only depends on $A_1$.
This bound follows from a computation that implies that the area of a minimal disk that is spanned by the given loop has to go to zero before time $B(A_1) \alpha_0 t$.
As a result, Hamilton's minimal disk argument yields a contradiction if $\alpha_0$ is chosen small enough such that $B(A_1) \alpha_0 t < t$.
This implies that the diameter of each $S_{i,t}$ must be bounded by $D_0 \sqrt{t}$ and hence the curvature is bounded by $K_0 t^{-1}$ on each $S_{i,t}$.
We hence obtain the desired curvature bound on the last remaining part of $M$, concluding the proof of the Main Theorem.
Note that the choice of constants $D_0, K_0$ and $\alpha_0$ is highly non-trivial.

\subsection{Structure of the following papers}
The proof of the two main Theorems \ref{Thm:MainTheorem-III} and \ref{Thm:geombehavior} is divided into the following four papers:

\begin{paperA}
In this paper we define precisely what we mean by ``Ricci flows with surgery and precise cutoff''.
The definition is chosen such that it incorporates most of the common notions of Ricci flows with surgery.
Then we review Perelman's analysis of these flows.
We will carry out some of Perelman's arguments again and generalize them to the case in which the underlying manifold is non-compact or has a boundary.
Then we prove various generalizations of Perelman's long-time estimates in the collapsing case.
These include several versions and localizations of the Generalized Key Lemma as mentioned above.
\end{paperA}

\begin{paperB}
This paper deals with area bounds of minimal surfaces or simplicial complexes under the Ricci flow.
We will show that minimal simplicial complexes are bounded by $A t$ in area and we will recall Hamilton's minimal disk argument, as well as a version of his argument for spheres.

The methods used in this paper are purely analytical.
Surgeries only play a very minor role.
\end{paperB}

\begin{paperC}
In this paper, we first recall several facts from the topology of $3$-manifolds.
Then we construct the (topological) simplicial complex to which we will apply the area estimate from the previous paper.
The simplicial complex is chosen in such a way that, in many cases, it is possible to extract a disk from its intersection with an arbitrary incompressible solid torus.
This fact is surprisingly non-trivial and it's proof occupies the major part of the paper.

Note that this paper is purely topological in nature.
Ricci flows will not be used.
\end{paperC}

\begin{paperD}
In this paper we assemble the results obtained in the previous papers and we finish the proofs of Theorems \ref{Thm:MainTheorem-III} and \ref{Thm:geombehavior}.
The paper contains the geometric characterization of the thin part $M_{\textnormal{thin}} (t)$ and a topological discussion of its decomposition.
We will then understand its fragmentation into good and bad parts and use the results of paper A to deduce the desired curvature bound first outside the solid tori $S_{1,t}, \ldots, S_{m_t,t}$, as mentioned in the outline, then on those $S_{i,t}$ of bounded diameter and eventually on collar neighborhoods of those $S_{i,t}$ of large diameter.
Finally, we use the simplicial complex from paper C together with the area bound of paper B to construct the required disk of bounded area.
This will then yield the desired contradiction.
\end{paperD}

We mention that each of the papers A, B and C is essentially self-explanatory and doesn't use any results of the other papers.
Only paper D makes use of the results of papers A--C.

\subsection{Acknowledgments}
I would like to thank Gang Tian for his constant help and encouragement, John Lott for many long conversations and Richard Schoen for pointing out Choe's work to me.
I am also indebted to Bernhard Leeb and Hans-Joachim Hein, who both fundamentally contributed to my understanding of Perelman's work.
Thanks also go to Simon Brendle, Alessandro Carlotto, Will Cavendish, Otis Chodosh, Daniel Faessler, Robert Kremser, Tobias Marxen, Rafe Mazzeo, Hyam Rubinstein, Andrew Sanders, Stephan Stadler and Brian White.

\end{document}